\documentclass[11pt]{article}
\usepackage{amsmath, amscd, amsfonts, amssymb, amsthm}

\newtheorem{theorem}{Theorem}[section]

\newtheorem{lemma}[theorem]{Lemma}
\newtheorem{prop}{Proposition}[section]

\def\fns{\scriptsize}
\def\mn{\medskip\noindent}

\begin{document}

\title{Two Phase Transitions for the \\
Contact Process on Small Worlds}
\author{Rick Durrett and Paul Jung}

\maketitle

\begin{abstract}
In our version of Watts and Strogatz's small world model, space is
a $d$-dimensional torus in which each individual has in addition
exactly one long-range neighbor chosen at random from the grid.
This modification is natural if one thinks of a town where an
individual's interactions at school, at work, or in social
situations introduce long-range connections. However, this change
dramatically alters the behavior of the contact process, producing
two phase transitions. We establish this by relating the small
world to an infinite ``big world" graph where the contact process
behavior is similar to the contact process on a tree.
\end{abstract}


\section{Introduction}
Small world graphs were first introduced by Watts and Strogatz
(1998).  In their model, they take a one-dimensional ring lattice
and connect all pairs of vertices that are distance $m$ or less.
They then ``rewire" each edge with probability $p$ by moving one
of the ends at random, where the new end is chosen uniformly. This
leads to a graph that has small diameter but, in contrast to the
Erd\"os-Renyi model, has a nontrivial density of triangles. These
are both properties that they observed in the collaboration graph
of film actors, the power grid, and the neural network of the
nematode, {\it C. elegans.}

The small world model has been extensively studied, although most
investigators have found it more convenient to study the Newman
and Watts (1999) version in which short-range connections are not
removed and each neighbor has a connection to a long-range
neighbor with probability $p$ (See Figure 1). Its graph theoretic
properties (e.g., the average distance between two points and the
clustering coefficient) are well understood (see Albert and
Barab\'asi (2002) for the physicist's view point or Barbour and
Reinert (2001) for rigorous results). Our focus here will be on
the behavior of processes taking place on these networks. Chapter
six of Watts(1999) discusses the SIR
(susceptible-infected-removed) disease model on small world graphs
in which individuals that are {\it susceptible} (state 0) become
{\it infected} (state 1) at a rate proportional to the number of
infected neighbors. Infected individuals, after a random amount of
time of fixed distribution, become {\it removed} (state 2), i.e.,
forever immune to further infection.

The SIR model on the small world graph has a detailed theory due
to its connection to percolation: we draw an oriented edge from
$x$ to $y$ if $x$ will succeed in infecting $y$ during the time it
is infected and the persistence of the epidemic is equivalent to
percolation.  See Section 8.2 of Newman (2003) for what is known
about the SIR models on small world graphs. Here we will
investigate the more difficult SIS
(susceptible-infected-susceptible) model, known to probabilists as
the contact process, where recovered individuals are immediately
susceptible. Moreno, Pastor-Satorras and Vespigiani (2002) have
studied this model by simulation, but we know of no rigorous
results for the contact process on the small world. Berger, Borgs,
Chayes, and Saberi (2004) have proved rigorous results for the
contact process on the preferential attachment graph of Barab\'asi
and Albert (1999).

Our version of the small world will be as follows.  We start with
a $d$-dimensional torus $(\mathbf{Z} \bmod R)^d$ and connect all
pairs of vertices within distance $m$ of each other using the
$\|\cdot\|_\infty$ norm. We require $R$ to be even so that we can
partition the $R^d$ vertices into $R^d/2$ pairs. Consider all such
partitions and then pick one at random. A new edge is then drawn
between each pair of vertices in the chosen partition. When $m\ll
R$, we think of these new edges as long-range connections. We will
call this graph $\mathcal{S}_m^R$, keeping in mind that this is a
random graph.

The reason for insisting that all individuals have exactly one
long-range neighbor is that we can define an associated ``big
world" graph $\mathcal{B}_m$ that is non-random. We define $\cal{
B}_m$ to consist of all vectors $\pm (z_1,\ldots,z_n)$ with $n\ge
1$ components with $z_j \in \mathbf{Z}^d$ and $z_j \neq 0$ for
$j<n$. Neighbors in the positive half-space are defined as
follows: a point $+(z_1,\ldots,z_n)$ is adjacent to
$+(z_1,\ldots,z_n+y)$ for all $y$ with $0 < \|y\|_\infty \le m$
(these are the short-range neighbors of $+(z_1,\ldots,z_n)$).  The
long-range neighbor is
\begin{eqnarray*}
+(z_1,\ldots,z_n,0)&&\text{ if }z_n\neq 0\\
+(z_1,\ldots,z_{n-1})&&\text{ if }z_n=0, n>1\\
-(0)\qquad &&\text{ if }z_n=0, n=1.
\end{eqnarray*}
See Figure 2 for a picture of the one dimensional case with $m=1$.
Of course in this case the graph is a tree.

The drawing of the small world in Figure 2 is more convenient for
our proof, but the graph can be more succinctly described if all
the long-range edges point down. In this case it is the free
product $\mathbf{Z}^d * \{0,1\}$ where the second factor is
$\mathbf{Z}$ mod 2. Elements of this group have the form $z_0 \, 1
\, z_1 \, 1 \, \ldots 1 \, z_n$ where $z_i \in \mathbf{Z}^d-\{0\}$
for $0 < i < n$. In words, this is the point you reach by moving
$z_0$ in the first copy, going down a long-range edge, moving
sideways by $z_1$, going down a long-range edge, etc.

We will consider the discrete-time contact process. On either the
small world or the big world, an infected individual lives for one
unit of time. During its infected period it will infect some of
its neighbors. All infection events are independent, and each site
that receives at least one infection is occupied with an infected
individual at the next time. A site infects itself or its
short-range neighbors, each with probability $\alpha/(2m+1)^d$. It
infects its long-range neighbor with probability $\beta$. Let
$\lambda=\alpha+\beta$ and $r=\alpha/\beta$. Hereafter we will
assume $\alpha>\beta$, and we fix the ratio $1<r<\infty$. We will
use $B_t$ to denote the contact process on the big world and
$\xi_t$ for the contact process on the small world.

The number of sites within distance $n$ of a given site in the big world
grows exponentially fast, so it is natural to guess that its
contact process will behave like the contact process on a
tree. Consider a tree in which each vertex has the same degree,
let $0$ be a distinguished vertex (the origin) of
the tree, and let $A^0_t$ be the set of infected sites at
time $t$ on the tree starting from 0 occupied. For the contact
process on the tree or on the big world, we can define two
critical values:
\begin{eqnarray}\label{sm2}
\lambda_1=\inf\{\lambda:\mathbf{P}(|A_t^0|=0
\text{ eventually})<1\}\\
\lambda_2=\inf\{\lambda:\liminf_{t\rightarrow\infty}\mathbf{P}(0\in
A_t^0)>0\}\nonumber.
\end{eqnarray}
We call $\lambda_1$ the weak survival critical value and
$\lambda_2$ the strong survival critical value. Pemantle (1992)
showed that for homogeneous trees where every vertex has at least
four neighbors, $\lambda_1 < \lambda_2$. He and Liggett (1996) who
extended the result to trees with degree 3, did this by finding
numerical bounds on the two critical values which showed they were
different. Later Stacey (1996) found a proof that did not rely on
numerical bounds. For another approach to separating the critical
values see Lalley and Selke (2002).

\begin{theorem}\label{theorem}
For each $0<r<\infty$ there exists $m_0$ such that for all $m\ge
m_0$, $\lambda_1<\lambda_2$ for the contact process on
$\mathcal{B}_m$.
\end{theorem}

The proof of Theorem \ref{theorem} will follow from Propositions
\ref{prop1} and \ref{prop3} below. Our approach to proving Theorem
\ref{theorem} is to find an upper bound $U$ (Proposition
\ref{prop1}) on the limiting weak survival critical value
$\limsup_{m\rightarrow\infty}\lambda_1$, a lower bound $L$
(Proposition \ref{prop3}) on the limiting strong survival critical
value $\liminf_{m\rightarrow\infty}\lambda_2$, and then to show
that $U<L$.

To obtain the two bounds, we will compare the contact process to
the branching random walk which has births like the contact
process, but allows more than one particle at a site. More
explicitly, each individual lives for one unit of time and sends
offspring to its site and to its short-range neighbors with
probability $\alpha/(2m+1)^d$ and to long-range neighbors with
probability $\beta$.  Birth events are independent, and the state
of a given site on the next time step is the total number of
births there.

The number of particles of the branching random walk at a given
site stochastically dominates the number of particles of the
contact process at a given site.  Therefore weak (strong) survival
of the contact process on $\mathcal{B}_m$ implies weak (strong)
survival of the branching random walk on $\mathcal{B}_m$. However,
in order to get an upper bound on the weak survival critical value
of the contact process, we need the opposite implication.  This is
achieved by showing that in the limit as $m\to\infty$, the
behavior of the contact process is essentially the same as that of
the branching random walk.  The following proposition, proved in
Section \ref{ub}, shows that the trivial necessary condition for
weak survival gives the asymptotic boundary of the survival
regime.

\begin{prop}If $\alpha+\beta> 1$ then the contact process survives for large
$m$. \label{prop1}
\end{prop}

To obtain a lower bound $L$ on $\lambda_2$, we use the fact that
strong survival of the contact process on $\mathcal{B}_m$ implies
strong survival of the branching random walk on $\mathcal{B}_m$.
Let $\lambda_2^{brw}(m)$ be the strong survival critical value of
the branching random walk. To compute the limit of
$\lambda_2^{brw}(m)$, we define the ``comb" of degree $m$,
$\mathcal{C}_m$, by restricting $\mathcal{B}_m$ to vertices of the
form $\{+(z),+(z,0),-(0)\}$ and all edges between any of these
vertices. As before, $+(z)$ and $+(z,0)$ are long-range neighbors
as are $+(0)$ and $-(0)$. The short-range neighbors of $+(z)$ are
$+(z+y)$ for $0 < \|y\|_\infty \le m$. The vertices $+(z,0)$ and
$-(0)$ have no short-range neighbors. To see the reason for the
name look at Figure 3 which gives a picture of the graph for $m=1,
d=1$. As we will show in Section \ref{sectionlb},
$\lim_{m\rightarrow\infty}\lambda_2^{brw}(m)$ is the critical
value for survival of the branching random walk on the comb.

It will be shown in Section \ref{sectionlb} that the branching
random walk on the comb survives if the largest eigenvalue of
$$
\begin{pmatrix} \alpha & \beta \cr \beta & 0 \end{pmatrix}
$$
is larger than 1 (since the comb is like $\mathbb{Z}^d$ the weak
survival and strong survival critical values are the same).
Solving the quadratic equation
$(\alpha-\lambda)(-\lambda)-\beta^2=0$ the largest root is
$$
\frac{ \alpha+ \sqrt{\alpha^2 - 4 \beta^2} }{2}
$$
A little algebra shows that this is larger than $1$ exactly when
$\alpha + \beta^2 > 1$. The proof of the following proposition is
given in Section \ref{sectionlb}. Since the critical value for the
branching random walk is smaller on the comb than it is on the big
world, this result is sharp asymptotically.

\begin{prop}\label{prop3} If $\alpha + \beta^2 < 1$ then there is no
strong survival in the contact process for large $m$.
\end{prop}

\mn  Comparing the above with Proposition \ref{prop1} shows that
we have separated the critical values for all $\beta<1$.

We have not been able to generalize Stacey's elegant argument to
the big world graph. However, we have been able to establish the
following ingredient used in his argument, which was proved for
the tree by Morrow, Schinazi, and Zhang (1994). Here, and in what
follows, we will write 0 as shorthand for $+(0)$.

\begin{theorem}\label{theorem2}
There exist constants $C_1>0$ and $C_2$ such that
$$
\exp(C_2 t)\le\mathbf{E}(|B^0_t|)\le \frac{1}{C_1}\exp(C_2 t).
$$
Moreover, $C_2<0$ when $\lambda<\lambda_1$, $C_2=0$ when
$\lambda=\lambda_1$, and $C_2>0$ when $\lambda>\lambda_1$.
\end{theorem}

\noindent Let $\tau_{\mathcal{B}}=\min\{t: B_t^{0}=\emptyset\}$ be
the extinction time of the contact process on the big world. It
follows from Theorem \ref{theorem2} that
$\tau_{\mathcal{B}}<\infty$ with probability one when $\lambda\le
\lambda_1$ and when $\lambda< \lambda_1$, $\tau_{\mathcal{B}}$ has a geometric
tail.

Having established the existence of two phase transitions on the
big world, our next question is: How does this translate into
behavior of the contact process on the small world? Let
$\sigma_{\mathcal{B}}=\min\{t: B_t^{0}=\emptyset\text{ or }0\in
B_t^0 \}$ be the first time that the infection either dies out or
comes back to the origin starting from one infection there at time
$0$. When $\lambda_1<\lambda<\lambda_2$, $\tau_\mathcal{B}$ and
$\sigma_\mathcal{B}$ are both infinite with positive probability,
and when $\lambda>\lambda_2$, $\sigma_\mathcal{B}$ is
almost-surely finite. Let $\tau_{\mathcal{S}}=\min\{t:
\xi_t^{0}=\emptyset\}$ and $\sigma_{\mathcal{S}}=\min\{t\ge 1:
\xi_t^{0}=\emptyset\text{ or }0\in \xi_t^0 \}$ be the
corresponding times for the contact process on the small world.

\begin{theorem}\label{cor}
Writing $\Rightarrow$ for convergence in distribution as $R \to
\infty$ we have \break $(a)$ $\tau_{\mathcal{S}}$ is
stochastically bounded above by $\tau_\mathcal{B}$ and
$\tau_{\mathcal{S}}\Rightarrow\tau_\mathcal{B}. \hfil$ \break
$(b)$ $\sigma_{\mathcal{S}}$ is stochastically bounded above by
$\sigma_\mathcal{B}$ and
$\sigma_{\mathcal{S}}\Rightarrow\sigma_\mathcal{B}$.
\end{theorem}

Theorems \ref{theorem} and \ref{cor} show that the small world
contact process exhibits different behavior in the three regimes
$(0,\lambda_1)$, $(\lambda_1,\lambda_2)$, and
$(\lambda_2,\infty)$. In the first interval the process dies out.
In the second, $\sigma_{\cal B} = \infty$ with positive
probability, so the time for the infection to return to 0 is not
tight (in the sense of weak convergence). In the third interval,
$\sigma_{\cal B} < \infty$ with probability one, so the infection
returns to 0 in a time that has a limit as the system size
increases. This is especially interesting since other researchers
have studied the contact process on the small world by simulation
without having noticed this qualitative difference between
$(\lambda_1,\lambda_2)$, and $(\lambda_2,\infty)$.

Since the small world is a finite graph, the infection will
eventually die out. However, by analogy with results for the
$d$-dimensional contact process on a finite set, we expect that if
the process does not become extinct quickly, it will survive for a
long time. Durrett and Liu (1988) showed that the supercritical
contact process on $[0,R)$ survives for an amount of time of order
$\exp(cR)$ starting from all ones, while Mountford (1999) showed
that the supercritical contact process on $[0,R)^d$ survives for
an amount of time of order $\exp(cR^d)$. At the moment we are only
able to prove the last conclusion for the following modification
of the small world contact process: each infected site infects its
short-range neighbors with probability $\alpha/(2m+1)^d$ and its
long-range neighbor with probability $\beta$, but now in addition,
it infects a random neighbor (uniformly) from the grid with
probability $\gamma>0$.

From a modelling point of view, this mechanism is reasonable. In
addition to long-range connections with friends at school or work,
one has random encounters with people one sits next to on
airplanes or meets while shopping in stores. In the language of
physics, the model with $\gamma=0$ has a quenched (i.e., fixed)
random environment, while the model with $\beta=0$ has an annealed
environment.

Our strategy for establishing prolonged survival is to show that
if the number of infected sites drops below $\eta R^d$, it will
with high probability rise to $2\eta R^d$ before dying out. To do
this we use the random connections to spread the particles out so
that they can grow independently. Ideally we would use the
long-range connections (instead of the random connections) to
achieve this; however, we have to deal with unlikely but annoying
scenarios such as all infected individuals being long-range
neighbors of sites that are respectively short-range neighbors of
each other.

\begin{theorem}\label{theorem1.4}
Consider the modified small world model on a $d$-dimensional torus
of radius $R$ with random infections at rate $\gamma>0$. If
$\lambda > \lambda_1$ and we start with all infected individuals
then there is a constant $c>0$ so that the probability the
infection persists to time $\exp(c R^d)$ tends to 1 as
$R\to\infty$.
\end{theorem}

\noindent This result is somewhat surprising and is qualitatively
different from the one Stacey (2001) proved for a finite
homogeneous tree with radius $R$. In that case if one starts with
all sites occupied, the extinction time grows linearly in $R$ for
$\lambda < \lambda_2$, but grows doubly exponential in $R$ for
$\lambda > \lambda_2$. In contrast Theorem \ref{theorem1.4} is, to
our knowledge,  the first result on any graph for which the
survival time of the contact process is exponential in volume when
$\lambda_1<\lambda<\lambda_2$. Of course the reason we can get
such a result is because a truncated tree is much different from
our locally tree-like small world.

We should mention that Berger, Borgs, Chayes, and Saberi (2004)
have studied the contact process on Barab\'asi and Albert's (1999)
scale-free graph, which has a power law degree distribution $\sim
Ck^{-3}$ and a critical value $\lambda_c=0$. For all fixed
$\lambda > 0$, they were able to show survival up to time
$\exp(c\lambda^2 n^{1/2})$ for the contact process on a graph with
$n$ vertices, but again, this is sub-exponential in volume. To do
this they combined an estimate for the survival time of the
contact process on a star graph with large degree, with the fact
that the largest degree in their random graph is $O(n^{1/2})$.

\medskip
{\bf Acknowledgments.} The authors would like to thank a referee
for suggesting that the results be formulated in terms of
the two parameters $\alpha$ and $\beta$, an associate
editor for indicating the connection between the big
world and the free product, and Laurent Saloff-Coste for
a tutorial on free products.

\section{An Upper Bound on Weak Survival}\label{ub}
To prove Proposition 1.1, consider the branching random walk $B_t$
on the half-space $\mathcal{B}_m^+$ consisting of all points of
the form $+(z_1,\ldots,z_n)$. We say a particle is at level $l$ if
it is at a point of the form $+(z_1,\ldots,z_l)$.

Clearly, the branching random walk is supercritical when
$\alpha+\beta>1$. The idea is to run the supercritical process to
time $t$ large enough so that we get enough special particles at
different points all of the form $+(z_1,\ldots,z_{l-1},0)$ on some
level $l$. We can then consider the half-spaces formed by all
points at level greater than or equal to $l$ treating points
$+(z_1,\ldots,z_{l-1},0)$ as the origins. Each of the special
particles starts its own branching random walk on these disjoint
half-spaces, and by construction none of these branching random
walks will ever run into each other (although a given branching
random walk can run into itself). If the expected number of
special particles at time $t$ is greater than 1, then a comparison
with a branching process shows that the total number of special
particles on levels $nl, n\in\mathbb{N}$ at times $nt$ is infinite
with positive probability. We then choose $m$ large so that the
contact process looks like a branching random walk.

Let us now carry out this plan more rigorously. Letting
$$\tilde{\alpha}=\alpha\left(1-\frac{1}{(2m+1)^d}\right),$$ we have that $\mathbf{E}(|B_t|) \ge
\tilde{\alpha}(\tilde{\alpha}+\beta)^{t-1}$.  One way to see this
is by disallowing short-range births onto points of the form,
$+(z_1,\ldots,z_{\ell(T)-1},0)$. The particles in $B_t$ are on
levels $0, \ldots, t$ so there is a level $\ell(t)$ where the
expected number of particles is at least
$\tilde{\alpha}(\tilde{\alpha}+\beta)^{t-1}/(t+1)$. When
$\alpha+\beta>1$ we can choose $m$ and $T$ large enough so that
\begin{eqnarray*}
\frac{\tilde{\alpha}(\tilde{\alpha}+\beta)^{T-1}}{T+1}
> \frac{3}{\tilde{\alpha}\beta}.\label{alphatilde}
\end{eqnarray*}

Modify the dynamics so that the only births occurring at time
$T+1$ are those across short-range edges on level $\ell(T)$
(excluding births to points of the form
$+(z_1,\ldots,z_{\ell(T)-1},0)$), and the only births occurring at
time $T+2$ are those across long-range neighbors from level
$\ell(T)$ to level $\ell(T)+1$. These birth at time $T+2$ will
lead to special particles on level $\ell(T)+1$, and the
half-spaces obtained by eliminating the long-range edges just
crossed are all disjoint.

The expected number of special particles in the modified branching
random walk on level $\ell(T)+1$ at time $T+2$ is at least 3. We
claim that if $m$ is large then the expected number of special
particles on level $\ell(T)+1$ at time $T+2$ in the contact
process on the half space is at least 2.

The expected number of particles in the branching process at any
time $t\le T$ is at most $({\alpha}+\beta)^T$. Let $E_K$ be the
event that $|B_t|\le K(\alpha+\beta)^T$ for all $t\le T$. For
$\epsilon>0$ choose $K$ large enough so that
$$\mathbf{E}(|B_T|1_{E_K})>(1-\epsilon)\mathbf{E}(|B_T|).$$

We now compare the contact process to the branching random walk.
To go from time $t$ to time $t+1$ in the contact process, we first
let births occur across long-range edges. This cannot lead to two
generation $t+1$ particles landing on the same site. We then let
the short-range births occur one at a time. On the event $E_K$,
the probability of a collision on any given birth is at most
$K({\alpha}+\beta)^T/(2m+1)^d$, so the probability of a collision
on some birth at some time $t\le T$ is at most
$TK^2(\alpha+\beta)^{2T}/(2m+1)^d$.

Now taking $m$ large enough, the expected number of special
particles at level $\ell(T)+1$ at time $T+2$ in the contact
process is at least 2. Since these particles sit atop half spaces
that are disjoint, the number of special particles on level
$2(\ell(T)+1)$ at time $2(T+2)$ for these processes in their
disjoint half-spaces are independent. By comparing with a
branching process we conclude that the contact process survives.

\section{A Lower Bound on Strong Survival}\label{sectionlb}

The first step is to compute the critical value for the branching
random walk on a comb. This is not necessary for the proof of the
lower bound on strong survival but it explains the expression for
the bound that we obtain. Additionally, it motivates us for the
latter part of this section. Here we will return to our original
notation and fix $r = \alpha/\beta
>1$ and vary $\lambda = \alpha+\beta$.

\begin{lemma}\label{combcr}
The critical value for the branching random walk on the comb
${\cal C}_m$ is independent of $m$ and is equal to
$$
\frac{2(r+1)}{r+\sqrt{r^2+4}}
$$
\end{lemma}

\begin{proof} We consider sites $+(z)$ to be level 1 and all other sites on the comb to be
level 2.  The level transition probabilities for a random walk on
the comb are given by
$$\left(
\begin{array}{cc}
\frac{\alpha}{\alpha+\beta}&\frac{\beta}{\alpha+\beta}\\
\frac{\beta}{\alpha+\beta}&0
\end{array}
\right)=\left(
\begin{array}{cc}
\frac{r}{r+1}&\frac{1}{r+1}\\
\frac{1}{r+1}&0
\end{array}
\right)
$$
Solving the quadratic equation
\begin{equation}
\label{quadeq} 0 = \left( \frac{r}{r+1} - x \right) (-x) -
\frac{1}{(r+1)^2} = x^2 - \frac{r}{r+1} x - \frac{1}{(r+1)^2}\  ,
\end{equation}
the largest eigenvalue of the matrix is
$$
\frac{r + \sqrt{r^2 + 4}}{2(r+1)}.
$$
This shows that the expected number
of particles at 0 tends to 0 exponentially fast if
$$
\lambda < \frac{2(r+1)}{r + \sqrt{r^2 + 4}}
$$
and gives the lower bound on the critical value. For the other
direction, we note that if the random walk is on level 1 after $n$
steps then there have been $m$ down steps, $m$ up steps, and
$n-2m$ sideways steps. The sideways steps are independent of the
up and down steps so the probability of being at 0 after $n$
steps, given that the walk is on level $0$ is greater than $
c/n^{d/2}$. This shows that the expected number of particles at 0
grows exponentially fast if
$$
\lambda > \frac{2(r+1)}{r + \sqrt{r^2 + 4}}
$$
The proof can now be completed using techniques of
Madras and Schinazi (1992) or Pemantle and Stacey (2001).
\end{proof}

Since the branching random walk $\eta^B_t$ on $\mathcal{B}_m$ dominates
the branching random walk $\eta^C_t$ on the comb of degree $m$, we have
$$
\lambda_2^{brw}(\mathcal{C}_m) \ge  \lambda_2^{brw}(\mathcal{B}_m)
$$
However, this inequality is in the opposite direction of what we
want. To get a result in the other direction, we introduce another
graph structure $\mathcal{K}_M$ which is similar to
$\mathcal{B}_m$ except that 0 has no short-range neighbors, and we
replace copies of $\mathbf{Z}^d$ by copies of the complete graph
on $M=(2m+1)^d$ vertices. Again, each point can be described
algebraically by a vector $(z_1,\ldots,z_n)$ with $n\ge 1$ integer
components, but now we have $z_1=0$ and $0\le z_j < M$ where
$z_j\neq 0$ for $1<j<n$. Figure 4 shows a picture of
$\mathcal{K}_M$ with $M=5$.

Consider two branching random walks $\eta^B_t$ and $\eta_t^K$ on
the respective graphs $\mathcal{B}_m$ and $\mathcal{K}_M$ starting
from one particle at the origin.  In the second process,
individuals that would be sent to short-range neighbors of 0 stay
at 0. We can couple $\eta^B_t$ and $\eta^K_t$ so that
$\eta_t^B(0)\le\eta_t^K(0)$ for all $t$.

Let $\lambda^{brw}_2(\mathcal{K}_M)$ denote the strong survival
critical value of $\eta^K_t$. In order to get a lower bound $L$
for the strong survival critical value in terms of $\lambda$ and
$r$, it is enough to show that
\begin{equation}\label{sm5}
\frac{2(r+1)}{r+\sqrt{r^2+4}}\le\liminf_{M\rightarrow\infty}\lambda^{brw}_2(\mathcal{K}_M).
\end{equation}
If we then interpret this in terms of $\alpha$ and $\beta$,
Proposition \ref{prop3} will have been proved.

Let $S_k$ be the random walk on $\mathcal{K}_M$ which jumps to
short-range neighbors with probability
$\alpha/((2m+1)^d(\alpha+\beta))$ and long-range neighbors with
probability $\beta/(\alpha+\beta)$. Jumps from 0 to short-range
neighbors make no change in the location of the particle. Define
$$
\tau=\inf\{k\ge 1:S_k=0\}\quad\text{and}\quad F^\lambda=\sum_{k\ge
1}\lambda^{k}\mathbf{P}(\tau=k).$$

\begin{lemma}\label{prop2}
$\lambda_2^{brw}(\mathcal{K}_M)\ge \sup\{\lambda:F^\lambda<1\}.$
\end{lemma}

\begin{proof}
For the branching random walk started with one particle at $0$:
$$
\mathbf{E}(\eta_k^K(0))=\lambda^k\mathbf{P}(S_k=0).
$$
Summing over $k$ we define
$$
G^{\lambda}=\sum_{k\ge 0}\lambda^k\mathbf{P}(S_k=0).
$$
We need only show that $F^\lambda<1$ implies $G^{\lambda}<\infty$, since this in turn
implies that the branching random walk dies out locally.
Breaking things down according to the value of $\tau$ we have for
$k\ge 1$
$$
\mathbf{P}(S_k=0)=\sum_{l=1}^k\mathbf{P}(\tau=l)\mathbf{P}(S_{k-l}=0).
$$
Multiplying by $\lambda^k$ and summing over $1\le k<\infty$ we get
$$
G^\lambda-1=F^\lambda G^\lambda
$$
so that $G^\lambda={1}/{(1-F^\lambda)}$.
\end{proof}

To examine the behavior of $S_k$ we will compare it to a birth and
death chain $R_k=\phi(S_k)$ on the nonnegative integers. $\phi$
maps 0 in $\mathcal{K}_M$ to 0, while for $n\ge 2$, a point of the
form $(z_1,\ldots,z_{n-1},0)$ maps to $2n-3$ and a point of the
form $(z_1,\ldots,z_n)$ with $z_n \neq 0$ maps to $2n-2$. See
Figure 4 for a picture of the mapping. Let $u = 1/(1+r)$ and $M =
(2m+1)^d$. It is clear from the symmetries of $\mathcal{K}_m$ that
$R_k$ is a Markov chain with transition probabilities: $r(0,0) = 1
- u$ and $r(0,1) = u$. When $j$ is odd, the last coordinate is 0
so:
$$
r(j,j-1) = u, \quad r(j,j) = (1-u)/M, \quad r(j,j+1) = (1-u)(1- 1/M).
$$
When $j\ge 2$ is even the last coordinate is non-zero so:
$$
r(j,j+1) = u, \quad r(j,j-1) = (1-u)/M, \quad r(j,j) = (1-u)(1- 1/M).
$$
It follows from our definitions that $\tau=\inf\{k\ge 1:R_k=0\}$.
For $x>0$, let
$$
h(x) = \sum_{k\ge 1}\lambda^{k}\mathbf{P}_x(\tau=k)
$$
and define $h(0)=1$. By considering what happens on one step, if
$x>0$
\begin{equation}
\label{hr} h(x) = \lambda \sum_y r(x,y) h(y).
\end{equation}
Using this equation for $x=2n, 2n-1$ and $2n+1$, we have for $n
\ge 1$
\begin{eqnarray*}
h(2n)&& \kern-2.75em \left(1-\frac{\lambda(1-u)}{M}\right)
 =  \lambda^2 h(2n-2) r(2n,2n-1) r(2n-1,2n-2) \\
& + & \lambda^2 h(2n+2) r(2n,2n+1) r(2n+1,2n+2) \\
& + & h(2n) \Bigl[ \lambda r(2n,2n)\left(1-\frac{\lambda(1-u)}{M}\right) + \lambda^2 r(2n,2n-1) r(2n-1,2n) \\
&& \hphantom{xxxxxxxxxxx} + \lambda^2 r(2n,2n+1) r(2n+1,2n)\Bigr].
\end{eqnarray*}
This can be rewritten as $ah(2n+2) - b h(2n) + c h(2n-2) = 0$
where
\begin{eqnarray*}
a & = & \lambda^2 u(1-u)(1-\frac{1}{M}) \\
b & = & 1-\frac{\lambda(1-u)}{M} - (\lambda-\frac{\lambda^2(1-u)}{M})(1-u)(1-\frac{1}{M})\\
&& \hphantom{xxxx} - \frac{\lambda^2 (1-u)^2(1-\frac{1}{M})}{M} - \lambda^2 u^2 \\
c & = & \frac{\lambda^2 u(1-u)}{M}.
\end{eqnarray*}
The solutions to the homogeneous difference equation are of the
form $h(2n) = C_1 \theta_1^n + C_2 \theta^n_2$ where $\theta_1$
and $\theta_2$ are the roots of $a\theta^2 - b \theta + c  = 0$,
or
$$
\theta_1 = \frac{ b + \sqrt{b^2-4ac} }{2a} \quad\hbox{and}\quad
\theta_2 = \frac{ b - \sqrt{b^2-4ac} }{2a}.
$$
When $\lambda<1$, $h(2n)$ is decreasing so $h(2n) = \theta_2^n$.
Since $h(x)$ is an analytic function of $\lambda$, this is the
formula for all $\lambda$ inside the radius of convergence. When
$M\to\infty$, $c \to 0$ while $a$ and $b$ have positive limits so
using the Maclaurin expansion of $\sqrt{1-x}$ we get
$$
\theta_2 = b \cdot \frac{ 1 - \sqrt{1 - \frac{4ac}{b^2} } }{2a}
\sim \frac{c}{b} \to 0
$$
and it follows that $h(2)\to 0$.

Intuitively what we have shown is that particles that reach level
2 (i.e., fall off the comb) and their descendants can be ignored.
To complete the calculation now we observe that using (\ref{hr})
with $x=1$ and letting $M\to\infty$ gives $h(1) \to \lambda u$. By
considering what happens on the first step starting from 0 we have
$$
F^\lambda = \lambda( 1-u + u h(1)) \to \lambda(1-u) + u^2
\lambda^2.
$$
Setting $F^\lambda=1$ and recalling $u = 1/(1+r)$ gives the quadratic equation
$$
\frac{1}{(1+r)^2} \lambda^2 + \lambda \frac{r}{r+1} - 1 = 0.
$$
The change of variable $\lambda=1/x$ gives (\ref{quadeq}) and
establishes (\ref{sm5}).

\section{Exponential Growth and Decay}
In this section we will prove Theorem \ref{theorem2}. Since we
have been writing $0$ for $+(0)$, we will similarly write $-0$ for
$-(0)$, and $-1$ for one of the nearest neighbors of $-(0)$. We
begin with a well-known property of trees, see e.g., Lemma 6.2 in
Pemantle (1992).

\begin{lemma}
\label{tree}Let $\mathbf{T}^3_r$ be the rooted binary tree in
which the root has degree one and all other vertices have degree
3. Suppose there are $I$ infected sites on $\mathbf{T}^3$. There
must be at least $I+2$ copies of $\mathbf{T}^3_r$, disjoint except
for possibly the root, where the only infected site is the root.
\end{lemma}

\begin{proof} We proceed by induction. If $I=1$ there are
$3=I+2$ copies. Each time we add a vertex we destroy one copy and
create two more.
\end{proof}

See Figure 5 for a picture of a case with $I=9$. We have drawn
$\mathcal{B}_1$ but it is the same as $\mathbf{T}^3$. As predicted
there are 11 copies of $\mathbf{T}^3_r$, marked by *'s on the
edges leading to their roots. We will call these edges, ``exterior
edges." The tree result extends to the following property of the
big world.

\begin{lemma}\label{bdry}
Given a finite set $B \subset \mathcal{B}_m$, there is a set of
vertices $D \subset B$ with $|D| \ge \frac{1}{c_d}|B|$ such that
each vertex $x\in D$ is either $(a)$ adjacent to a vacant
translate of a copy of $\mathcal{B}^+_m$ with $x$ sitting at $-0$,
or $(b)$ adjacent to a vacant translate of a copy of
$\mathcal{B}^+_m \cup \{-0\}$ with $x$ sitting at a site that is a
nearest neighbor of $-0$.
\end{lemma}

\begin{proof}
Consider first the case $d=1$. We can embed $\mathbf{T}^3$ into
$\mathcal{B}_m$. Using Lemma \ref{tree} it follows that there are
at least $(I+2)/3$ completely disjoint copies of $\mathbf{T}^3_r$
where the only infected site is the root. If the exterior edge
from the boundary vertex $x$ is a long-range edge then we are in
case (a). If the exterior edge from the boundary vertex $x$ is a
short-range edge we are in case (b). See Figure 5.

In $d=2$ we can still embed a tree but it has a variable degree.
In each plane we connect $(0,0)$ to its four nearest neighbors.
For $k\ge 1$ we connect $(k,0)$ to both $(k+1,0)$ and $(k,1)$
while we connect $(k,j)$ to only $(k,j+1)$ for $j\ge 1$. Then
extend the construction to the other four quadrants so that it is
symmetric under 90 degree rotations. See Figure 6 for a picture.
All points in addition have long distance neighbors so the result
holds with $c_2=5$. We leave the details for $d\ge 3$ to the
reader.
\end{proof}

With Lemma \ref{bdry} established, the rest of the proof is very
similar to the proof of Theorem 2 in Madras, Schinazi, and Zhang
(1994).

\begin{lemma}\label{madras}
For the contact process $B^0_t$ on $\mathcal{B}_m$, there exists
$C_1>0$ such that
$$\mathbf{E}(|B^0_{t+s}|)\ge C_1 \mathbf{E}(|B^0_{t}|)
\mathbf{E}(|B_{s}^0|).$$
\end{lemma}

\begin{proof}
For $j=0,1$, let $|B_t^{+,-j}|$ be the number of points in
$\mathcal{B}^+_m$ for the process started with a single infected
individual at $-j$. Let
$M_t=\min\{\mathbf{E}(|B_t^{+,-0}|),\mathbf{E}(|B_t^{+,-1}|)\}$.
By additivity, the Markov property, and Lemma \ref{bdry} we have
$$
\mathbf{E}(|B_{t+s}^0|) \ge \frac{1}{c_d}\mathbf{E}(|B_t^0|)M_s.
$$
The lemma will be proved when we show that there exists $C_1>0$
such that $M_s\ge C_1 \mathbf{E}(|B_s^0|)$.

To do this we begin by observing that for $j=0,1$
\begin{equation}
\label{extmean} \mathbf{E}(|B_t^{+,-j}|)\ge\frac{\beta}{2}\cdot
\frac{\alpha^2}{(2m+1)^{2d}} \cdot
\mathbf{E}(|B^{\{0,-0\}}_{t-2}|)
\end{equation}
where we are considering the event that on the first step $-j$
infects $-0$, and on the second step $-0$ infects itself and $0$.
The factor of $1/2$ comes from the fact that when starting from
$\{-0,0\}$, the expected number of points in the positive
half-space is one-half the number in the whole space. To take care
of the differences in starting configuration and time, we note
$$
\mathbf{E}(|B^0_t|) \le \mathbf{E}(|B^{\{0,-0\}}_t|) \le \lambda^2
\mathbf{E}(|B^{\{0,-0\}}_{t-2}|).
$$
From this it follows that for $j=0,1$
$$
\liminf_{t\rightarrow\infty}
\frac{\mathbf{E}(|B^{+,-j}_t|)}{\mathbf{E}(|B^0_t|)} \ge
\frac{\beta}{2} \cdot \left(\frac{\alpha}{(2m+1)^d
\lambda}\right)^2
$$
\end{proof}

\begin{proof}[Proof of Theorem \ref{theorem2}]
By the additivity of the contact process and the transitivity of
$\mathcal{B}_m$,
$$\mathbf{E}(|B^0_{t+s}|)\le \mathbf{E}(|B^0_{t}|) \mathbf{E}(|B_{s}^0|).$$
Subadditivity then implies that
\begin{equation}\label{sm9}
\lim_{t\rightarrow\infty}\frac{1}{t}\log
\mathbf{E}(|B^0_{t}|)=\inf_{t>0}\frac{1}{t}\log
\mathbf{E}(|B^0_t|)=C_2
\end{equation}
exists so that $\exp(C_2 t)\le \mathbf{E}(|B^0_t|).$ A similar
argument using superadditivity of $C_1 \mathbf{E}(|B^0_{t}|)$ and
Lemma \ref{madras} shows that
\begin{equation}\label{sm10}
\lim_{t\rightarrow\infty}\frac{1}{t}\log
C_1\mathbf{E}(|B^0_{t}|)=\sup_{t>0}\frac{1}{t}\log
C_1\mathbf{E}(|B^0_t|)=C_2
\end{equation}
so that $\frac{1}{C_1}\exp(C_2 t)\ge \mathbf{E}(|B^0_t|).$

We turn now to the sign of $C_2$ in the different regimes. When
$\lambda>\lambda_1$ we use the well-known fact that on the event
$$\{|B_t^0|\ge 1 \text{ for all }k\ge 0\},$$
$\lim_{t\rightarrow\infty}|B_t|=\infty$ almost surely. So by the
upper bound $\frac{1}{C_1}\exp(C_2 k)\ge \mathbf{E}(|B^0_t|)$, we
see that $C_2>0$ when $\lambda>\lambda_1$.

When $\lambda<\lambda_1$, the work of Aizenman and Barsky (1987)
gives
$$\sum_{t\ge 0} \mathbf{E}(|B_t^0|)<\infty$$ so that $C_2<0$ in
this regime.

For $\lambda=\lambda_1$ we use a continuity argument. Note that
for any fixed $t>0$, $$\lambda\rightarrow \frac{1}{t}\log
\mathbf{E}(|B^0_{t}|)$$ is a continuous function of $\lambda>0$.
By (\ref{sm9}) and (\ref{sm10}), $\lambda\rightarrow C_2$ is lower
and upper semi-continuous and thus continuous. Therefore when
$\lambda=\lambda_1$ it must be that $C_2=0$.
\end{proof}

\section{Proofs of Theorem \ref{cor} and \ref{theorem1.4}}

\begin{proof}[Proof of Theorem \ref{cor}]
We start by viewing $\mathcal{B}_m$ as a ``covering space" of
$\mathcal{S}_m^R$. Given a realization of the random graph
$\mathcal{S}_m^R$, we mark a distinguished vertex which we
identify with the origin of $\mathcal{B}_m$.  From there we
identify the long-range edge of the distinguished vertex with the
long-range edge of the origin in $\mathcal{B}_m$. Similarly, we
identify the short-range edges of the two graphs together in such
a way so that graph distances on the underlying structure given by
$\mathcal{B}_m$ are preserved. Continuing in this manner we can
identify each vertex in $\mathcal{B}_m$ with some vertex in our
realization of $\mathcal{S}_m^R$.

Using this identification, couple together the processes $B^{0}_t$
and $\xi_t^0$ for each realization of the random graph
$\mathcal{S}_m^R$. By additivity of the contact process, it is
clear that $|B^{0}_t|\ge |\xi^0_t|$ which shows that
$\tau_\mathcal{S}$ is bounded above by $\tau_\mathcal{B}$. Choose
$\epsilon>0$ and $N>0$. Using the graph identification above, let
$G_{R,N}$ be the event (on the random graph $\mathcal{S}_m^R$)
that there is no cycle of length $2N$ which contains $0$.  For
fixed $N$, we have that $\mathbf{P}(G_{R,N})\to 1$ as $R\to\infty$
(if this is not clear to the reader, an explicit argument is given
in the proof of \ref{theorem1.4} below). If there are no cycles of
length $2N$ containing $0$ in the graph $\mathcal{S}_m^R$ then no
two vertices within distance $N$ of $0$ are identified together in
$\mathcal{B}_m$.

Choose $R_0$ so that $R\ge R_0$ implies
$\mathbf{P}(G_{R,N})>1-{\epsilon}$. Using the coupling to identify
points of $\mathcal{B}_m$ with $\mathcal{S}_m^R$, we have for all
$R\ge R_0$,
$$
\mathbf{P}(B_k^0=\xi_k^0 \text{ for all }k\le N)>1-\epsilon
$$
which proves
$\lim_{R\rightarrow\infty}\tau_\mathcal{S}=\tau_\mathcal{B}$. The
proof of
$\lim_{R\rightarrow\infty}\sigma_\mathcal{S}=\sigma_\mathcal{B}$
is similar and is therefore omitted.
\end{proof}

\begin{proof}[Proof of Theorem \ref{theorem1.4}]
Let $M = (2m+1)^d$ be the number of neighbors and let
$$
\delta = \frac{\alpha\gamma}{8M}
$$
The contact process on the big world with $\gamma=0$ is
supercritical so we can pick $T$ so that $\mathbf{E}(|\xi^0_T|)
\ge 5/\delta$ and we can pick $K$ so that if $\bar \xi^0_T$ is the
contact process restricted to the ball of radius $K$ around 0 then
$\mathbf{E}(|\bar\xi^0_T|) \ge 4/\delta$. Let $N_K$ be the number
of points in the ball of radius $K$ around 0 on the big world.
Pick $\eta$ so that
$$
\eta\cdot\frac{\alpha \gamma N_{2K}}{4M} \le \frac{1}{5}
$$

It suffices to show that if the number of occupied sites drops
below $\eta R^d$ occupied sites then with probability greater than
$1 - \exp(-b R^d)$ for some $b>0$, it will return to having more
than $2\eta R^d$ occupied sites.  To see this we just choose $c<b$
giving us a high probability of having $\exp(c n^d)$ successful
recoveries. Since each recovery takes at least one unit of time,
the result follows.

The first step towards achieving the above is to show that the
first time it drops below $\eta R^d$ it does not fall too far. To
get a lower bound on how far it falls, we only look at births from
sites onto themselves. These are independent events with
probability $\alpha/M$. A standard large deviations result for the
binomial implies that if there are at least $\eta R^d$ particles
alive at time $t$ then the probability of less than $\eta R^d
\alpha/2M$ alive at time $t+1$ is less than $\exp(-b_1 R^d)$ for
some $b_1>0$. Thus the probability this occurs at some time before
$\exp(b_1 R^d/2)$ is less than $\exp(-b_1 R^d/2)$.

After the number of occupied sites falls below $\eta R^d$, our
next step is to randomize the locations by having one time step in
which we only allow births induced by the parameter $\gamma$. If
we start with $\eta R^d \alpha/2M$ particles, then after these
particles give birth onto their randomly chosen neighbors (with no
other births allowed) we will have, on average, $\eta R^d
\gamma\alpha/2M$ particles.  In fact, if $t+1$ is the time at
which the number of occupied sites falls below $\eta R^d$ (but not
too far), then using the large deviations estimate for the
binomial once more we have that the probability of there being
less than $\eta R^d\gamma\alpha/4M$ particles at time $t+2$ is
less than $\exp(-b_2 R^d)$ for some $b_2>0$ (we will take care of
double counting below). Having used the randomized births in this
step, we ignore them for the rest of the proof and use the process
with $\gamma=0$ which we have assumed is supercritical.

We say that a site $x$ is \textit{good} if the small world is
identical to the big world inside a ball of radius $K$ around $x$.
Fixing $x$, the probability of a self-intersection (i.e., the
probability that $x$ is not a good site) when we generate the ball
of radius $K$ around $x$ is less than $N_K^3/n$. To see this, we
grow the ball around $x$ starting from just $x$ and its
short-range neighbors. Adding the long-range neighbor of $x$, we
see that the probability it causes a self-intersection is bounded
above by $N_K/n$. But we also have to be worried about the
short-range neighbors of this long-range neighbor causing a
self-intersection, so we increase the bound generously to
$N_K^2/n$. Now since there are at most $N_K$ neighbors to add, the
bound $N_K^3/n$ follows; therefore, the probability that $x$ is a
bad site tends to 0 as $n \to \infty$.

Now fix a realization of the small world graph such that the
fraction of bad points does not exceed $1/20$. We want to find a
subset $G$ of the randomized births that lie on good sites and for
which the evolutions of $\bar{\xi}^x_t, x\in G$ are independent.
To do this we imagine that the successful randomized births occur
one at a time (successful with respect to $\gamma$). We accept the
first birth if it lands on a good site. We accept the second birth
if it lands on a good site and it does not fall within the ball of
radius $2K$ around the first birth, and we continue in this
manner.

When there are $j$ successful randomized births accepted, the
total number of sites that lie within the balls of radius $2K$
around them does not exceed $j N_{2K}$.  When
$$j\le \eta R^d\gamma\alpha/4M,$$
the probability of the next successful randomized birth landing in
the forbidden zone is less than $1/5$ by our choice of $\eta$
above. Since we are on a realization of $\mathcal{S}_m^R$ whose
fraction of bad points does not exceed $1/20$, the probability of
discarding a successful randomized birth conditioning on $j$ prior
acceptances is less than $1/4$. Thus the number of births
discarded out of the first $\eta R^d\gamma\alpha/(8m+4)$
successful randomized births is bounded by a binomial with success
probability $p=1/4$. Our large deviations estimate implies that
with probability less than $\exp(-b_3R^d)$, the number of
randomized births remaining, $|G|$, is larger than
$$\frac{\eta R^d\gamma\alpha}{8M}=\delta\eta R^d.$$

By our choice of $G$, the particles in $G$ can grow independently
up to a distance of $K$, or in other words, the processes
$\bar{\xi}_t^x, x\in G$ evolve independently. Since
$\mathbf{E}(|\bar{\xi}_T^0|)>4/\delta$, a final large large
deviations estimate for bounded i.i.d. random variables tells us
that the probability of ending up with fewer than $2\eta R^d$
particles at time $T$ is less than $\exp(-b_4 R^d)$.  The total of
our error probabilities is less than $\exp(-b R^d)$ which
completes the proof.
\end{proof}

\section*{References}

\mn Aizenman, M. and Barsky, D. J. (1987). Sharpness of the phase
transition in percolation models. {\it Comm. Math. Phys.}
\textbf{108}, 489-526.

\mn Albert, R., and Barab\'asi, A. L. (2002). Statistical
mechanics of complex networks. {\it Rev. Mod. Phys.} {\bf 74},
47--97

\mn Barbour, A.D., and Reinert, G. (2001). Small worlds. {\it
Random Structures and Algorithms.} {\bf 19}, 54--74

\mn Barab\'asi, A., and Albert, R. (1999) Emergence of scaling in
random networks. {\it Science.} {\bf 286}, 509--512

\mn Berger, N., Borgs, C., Chayes, J. T., and Saberi, A. (2004).
On the spread of viruses on the internet.

\mn Durrett, R. and Liu, X. (1988). The contact process on a
finite set. \emph{Ann. Probab.} {\bf 16}, 1158-1173.

\mn Lalley, Steven P. and Sellke, T. M. (2002). Anisotropic
contact processes on homogeneous trees. \emph{Stoch. Proc. Appl.}
{\bf 101}, 163--183.

\mn Liggett, T. M. (1996). Multiple transition points for the
contact process on the binary tree. \emph{Ann. Probab.} {\bf 24},
1675-1710.

\mn Madras, N. and Schinazi, R. B. (1992). Branching random walks
on trees. \emph{Stoch. Proc. Appl.} {\bf 42}, 255-267.

\mn Moreno, Y., Pastor-Satorras, R., and Vespigniani, A. (2002).
Epidemic outbreaks in complex heterogenenous networks. {\it Eur.
Phys. J., B.} {\bf 26}, 521--529

\mn Morrow, G. J., Schinazi, R. B., and Zhang, Y. (1994). The
critical contact process on a homogeneous tree. \emph{J. Appl.
Probab.} {\bf 31}, 250-255.

\mn Mountford, T.S. (1999). Existence of a constant for finite
system extinction. {\it J. Stat. Phys.} {\bf 96}, 1331--1341

\mn Newman, M. E. J., Jensen, I., and Ziff, R. M. (2002).
Percolation and epidemics in a two-dimensional small world.
\emph{Phys. Rev. E} {\bf 65}, paper 021904.

\mn Newman, M. E. J. (2003). The structure and function of complex
networks. {\it SIAM Review.} {\bf 45}, 167--256

\mn Pemantle, R. (1992). The contact process on trees. \emph{Ann.
Probab.} {\bf 20}, 2089-2116.

\mn Pemantle, R. and Stacey, A. M. (2001). The branching random
walk and contact process on the Galton-Watson and nonhomogeneous
trees. \emph{Ann. Probab.} {\bf 29}, 1563-1590.

\mn Stacey, A. M. (1996). The existence of an intermediate phase
for the contact process on trees. \emph{Ann. Probab.} {\bf 24},
1711-1726.

\mn Stacey, A. M. (2001). The contact process on finite
homogeneous trees. {\it Prob. Theor. Rel. Fields.} {\bf 121},
551-576

\mn Stacey, A. M. (2003). Branching random walks on
quasi-transitive graphs. \emph{Combin. Probab. Comput.} {\bf 12},
345-358.

\mn Watts, D. J. (1999). \emph{Small worlds: The dynamics of
networks between order and randomness}, Princeton University
Press, New Jersey.

\mn Watts, D. J. and Strogatz, S. H. (1998). Collective dynamics
of `small-world' networks. \emph{Nature} {\bf 393}, 440-442.

\clearpage

\begin{center}
\begin{picture}(300,300)
\put(270,150){$\bullet$}\put(269,159){$\bullet$}\put(268,168){$\bullet$}\put(266,178){$\bullet$}
\put(264,187){$\bullet$}\put(260,195){$\bullet$}\put(256,204){$\bullet$}\put(252,212){$\bullet$}
\put(247,220){$\bullet$}\put(241,227){$\bullet$}\put(234,234){$\bullet$}\put(227,241){$\bullet$}
\put(220,247){$\bullet$}\put(212,252){$\bullet$}\put(204,256){$\bullet$}\put(195,260){$\bullet$}
\put(187,264){$\bullet$}\put(178,266){$\bullet$}\put(168,268){$\bullet$}\put(159,269){$\bullet$}
\put(150,269){$\bullet$}\put(141,269){$\bullet$}\put(132,268){$\bullet$}\put(122,266){$\bullet$}
\put(113,264){$\bullet$}\put(105,260){$\bullet$}\put(96,256){$\bullet$}\put(88,252){$\bullet$}
\put(80,247){$\bullet$}\put(73,241){$\bullet$}\put(66,234){$\bullet$}\put(59,227){$\bullet$}
\put(53,220){$\bullet$}\put(48,212){$\bullet$}\put(44,204){$\bullet$}\put(40,195){$\bullet$}
\put(36,187){$\bullet$}\put(34,178){$\bullet$}\put(32,168){$\bullet$}\put(31,159){$\bullet$}
\put(31,150){$\bullet$}\put(31,141){$\bullet$}\put(32,132){$\bullet$}\put(34,122){$\bullet$}
\put(36,113){$\bullet$}\put(40,105){$\bullet$}\put(44,96){$\bullet$}\put(48,88){$\bullet$}
\put(53,80){$\bullet$}\put(59,73){$\bullet$}\put(66,66){$\bullet$}\put(73,59){$\bullet$}
\put(80,53){$\bullet$}\put(88,48){$\bullet$}\put(96,44){$\bullet$}\put(105,40){$\bullet$}
\put(113,36){$\bullet$}\put(122,34){$\bullet$}\put(132,32){$\bullet$}\put(141,31){$\bullet$}
\put(150,31){$\bullet$}\put(159,31){$\bullet$}\put(168,32){$\bullet$}\put(178,34){$\bullet$}
\put(187,36){$\bullet$}\put(195,40){$\bullet$}\put(204,44){$\bullet$}\put(212,48){$\bullet$}
\put(220,53){$\bullet$}\put(227,59){$\bullet$}\put(234,66){$\bullet$}\put(241,73){$\bullet$}
\put(247,80){$\bullet$}\put(252,88){$\bullet$}\put(256,96){$\bullet$}\put(260,105){$\bullet$}
\put(264,113){$\bullet$}\put(266,122){$\bullet$}\put(268,132){$\bullet$}\put(269,141){$\bullet$}
\put(34,170){\line(1,1){100}}
\put(258,98){\line(-1,2){87}}
\put(33,152){\line(1,0){240}}
\put(152,33){\line(1,3){72}}
\put(222,55){\line(-1,1){167}}
\put(268,125){\line(-3,-1){194}}

\end{picture}
\end{center}
\bigskip\noindent
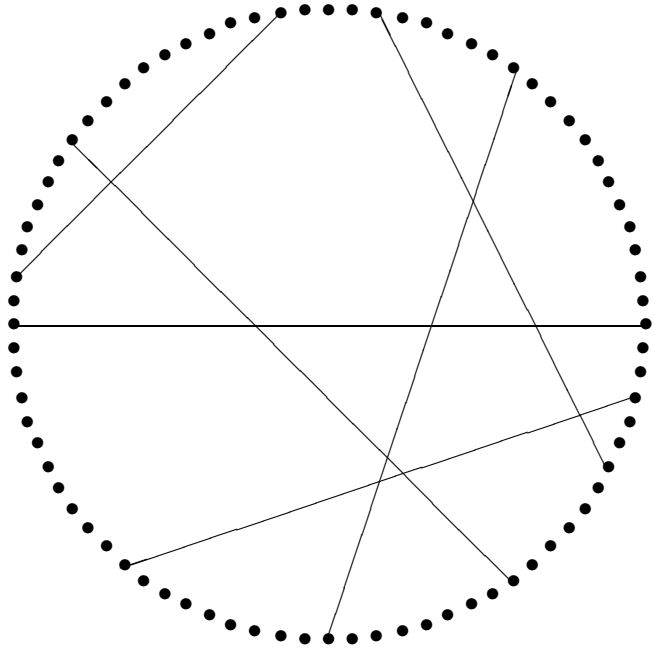
Figure 1. The Newman and Watts version of a small world graph. For
simplicity we have not drawn the short-range connections.

\clearpage

\begin{center}
\begin{picture}(300,270)
\put(150,110){\line(0,1){50}}
\put(25,110){\line(1,0){250}}
\put(25,160){\line(1,0){250}}
\put(50,160){\line(0,1){50}}
\put(100,160){\line(0,1){50}}
\put(200,160){\line(0,1){50}}
\put(250,160){\line(0,1){50}}
\put(50,110){\line(0,-1){50}}
\put(100,110){\line(0,-1){50}}
\put(200,110){\line(0,-1){50}}
\put(250,110){\line(0,-1){50}}
\put(25,185){\line(1,1){50}}
\put(75,185){\line(1,1){50}}
\put(175,185){\line(1,1){50}}
\put(225,185){\line(1,1){50}}
\put(25,35){\line(1,1){50}}
\put(75,35){\line(1,1){50}}
\put(175,35){\line(1,1){50}}
\put(225,35){\line(1,1){50}}
\put(34,194){\line(0,1){30}}\put(42,202){\line(0,1){30}}\put(58,218){\line(0,1){30}}\put(66,226){\line(0,1){30}}
\put(84,194){\line(0,1){30}} \put(92,202){\line(0,1){30}}\put(108,218){\line(0,1){30}}\put(116,226){\line(0,1){30}}
\put(184,194){\line(0,1){30}} \put(192,202){\line(0,1){30}} \put(208,218){\line(0,1){30}} \put(216,226){\line(0,1){30}}
\put(234,194){\line(0,1){30}} \put(242,202){\line(0,1){30}} \put(258,218){\line(0,1){30}} \put(266,226){\line(0,1){30}}
\put(34,44){\line(0,-1){30}} \put(42,52){\line(0,-1){30}} \put(58,68){\line(0,-1){30}} \put(66,76){\line(0,-1){30}}
\put(84,44){\line(0,-1){30}} \put(92,52){\line(0,-1){30}} \put(108,68){\line(0,-1){30}} \put(116,76){\line(0,-1){30}}
\put(184,44){\line(0,-1){30}} \put(192,52){\line(0,-1){30}} \put(208,68){\line(0,-1){30}} \put(216,76){\line(0,-1){30}}
\put(234,44){\line(0,-1){30}} \put(242,52){\line(0,-1){30}} \put(258,68){\line(0,-1){30}} \put(266,76){\line(0,-1){30}}
\put(140,100){\fns +(0)}
\put(90,117){\fns +(--1)}
\put(40,117){\fns +(--2)}
\put(190,117){\fns +(1)}
\put(240,117){\fns +(2)}
\put(140,165){\fns --(0)}
\put(90,150){\fns --(--1)}
\put(40,150){\fns --(--2)}
\put(190,150){\fns --(1)}
\put(240,150){\fns --(2)}
\put(20,62){\fns +(--2,0)}
\put(70,62){\fns +(--1,0)}
\put(189,78){\fns +(1,2)}
\put(179,69){\fns +(1,1)}
\put(170,60){\fns +(1,0)}
\put(159,51){\fns +(1,--1)}
\put(150,42){\fns +(1,--2)}
\put(222,62){\fns +(2,0)}
\put(267,44){\fns +(2,2,0)}
\put(259,36){\fns +(2,1,0)}
\put(243,20){\fns +(2,--1,0)}
\put(235,12){\fns +(2,--2,0)}
\end{picture}

\bigskip
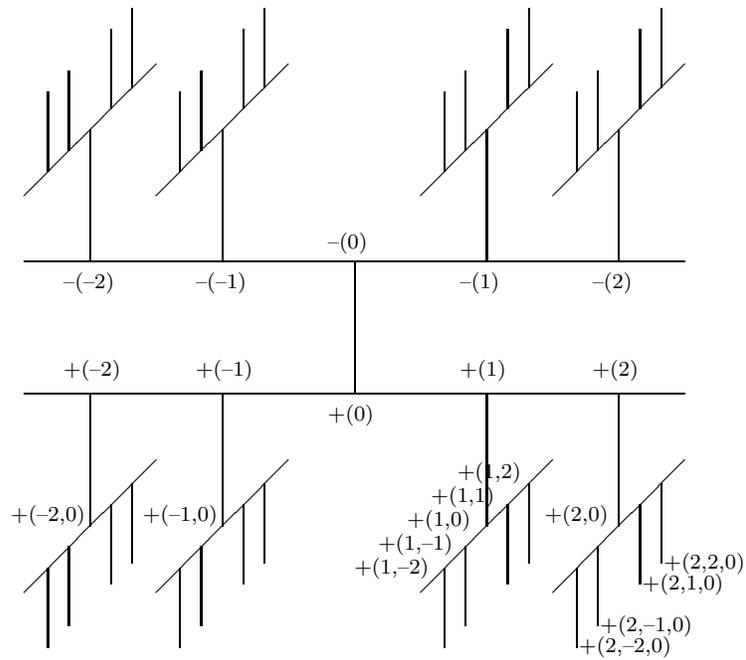
Figure 2. The big world graph, ${\cal B}_1$.
\end{center}

\clearpage

\begin{center}
\begin{picture}(310,100)
\put(30,80){\fns +(--3)}
\put(70,80){\fns +(--2)}
\put(110,80){\fns +(--1)}
\put(150,80){\fns +(0)}
\put(190,80){\fns +(1)}
\put(230,80){\fns +(2)}
\put(270,80){\fns +(3)}
\put(10,70){\line(1,0){290}}
\put(42,70){\line(0,-1){40}}
\put(82,70){\line(0,-1){40}}
\put(122,70){\line(0,-1){40}}
\put(162,70){\line(0,-1){40}}
\put(202,70){\line(0,-1){40}}
\put(242,70){\line(0,-1){40}}
\put(282,70){\line(0,-1){40}}
\put(30,20){\fns +(--3,0)}
\put(70,20){\fns +(--2,0)}
\put(110,20){\fns +(--1,0)}
\put(155,20){\fns --(0)}
\put(190,20){\fns +(1,0)}
\put(230,20){\fns +(2,0)}
\put(270,20){\fns +(3,0)}
\end{picture}
\end{center}

\bigskip\noindent
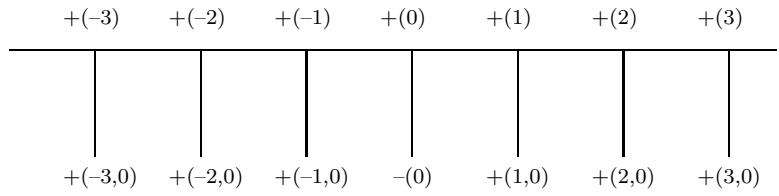
Figure 3. The comb $\mathcal{C}_1$. We have rotated the edge from $+(0)$ to $-(0)$ to make the picture better match the name.

\clearpage
\begin{center}
\begin{picture}(300,300)
\put(147,268){$\bullet$}
\put(156,268){0}
\put(150,270){\line(0,-1){50}}
\put(149,221){\line(-1,-1){21}}
\put(151,221){\line(1,-1){21}}
\put(134,177){\line(-1,4){6}}
\put(166,177){\line(1,4){6}}
\put(134,177){\line(1,0){33}}
\put(133,178){\line(2,5){17}}
\put(167,178){\line(-2,5){17}}
\put(129,200){\line(5,-3){36}}
\put(171,200){\line(-5,-3){36}}
\put(130,201){\line(1,0){40}}
\put(131,175){$\bullet$}
\put(164,175){$\bullet$}
\put(125,198){$\bullet$}
\put(169,198){$\bullet$}
\put(146,218){$\otimes$}
\put(156,223){1}
\put(177,203){2}
\put(117,203){2}
\put(121,175){2}
\put(174,175){2}
\put(133,178){\line(-2,-3){32}}
\put(167,178){\line(2,-3){32}}
\put(130,200){\line(-1,0){88}}
\put(170,200){\line(1,0){88}}
\put(99,130){\line(-1,-1){21}}
\put(101,130){\line(1,-1){21}}
\put(84,86){\line(-1,4){6}}
\put(116,86){\line(1,4){6}}
\put(84,86){\line(1,0){33}}
\put(83,87){\line(2,5){17}}
\put(115,87){\line(-2,5){17}}
\put(79,109){\line(5,-3){36}}
\put(121,109){\line(-5,-3){36}}
\put(80,110){\line(1,0){40}}
\put(81,84){$\bullet$}
\put(114,84){$\bullet$}
\put(75,107){$\bullet$}
\put(119,107){$\bullet$}
\put(95,127){$\otimes$}
\put(95,137){3}
\put(125,112){4}
\put(68,112){4}
\put(72,84){4}
\put(122,84){4}
\put(83,86){\line(-1,-1){12}}
\put(118,86){\line(1,-1){12}}
\put(78,109){\line(-1,0){16}}
\put(123,109){\line(1,0){16}}
\put(199,130){\line(-1,-1){21}}
\put(201,130){\line(1,-1){21}}
\put(184,86){\line(-1,4){6}}
\put(216,86){\line(1,4){6}}
\put(184,86){\line(1,0){33}}
\put(183,87){\line(2,5){17}}
\put(215,87){\line(-2,5){17}}
\put(179,109){\line(5,-3){36}}
\put(221,109){\line(-5,-3){36}}
\put(180,110){\line(1,0){40}}
\put(181,84){$\bullet$}
\put(214,84){$\bullet$}
\put(175,107){$\bullet$}
\put(219,107){$\bullet$}
\put(195,127){$\otimes$}
\put(197,137){3}
\put(225,112){4}
\put(168,112){4}
\put(172,84){4}
\put(222,84){4}
\put(183,86){\line(-1,-1){12}}
\put(218,86){\line(1,-1){12}}
\put(178,109){\line(-1,0){16}}
\put(223,109){\line(1,0){16}}
\put(39,201){\line(-1,-1){21}}
\put(41,201){\line(1,-1){21}}
\put(24,157){\line(-1,4){6}}
\put(56,157){\line(1,4){6}}
\put(24,157){\line(1,0){33}}
\put(23,158){\line(2,5){17}}
\put(57,158){\line(-2,5){17}}
\put(19,180){\line(5,-3){36}}
\put(61,180){\line(-5,-3){36}}
\put(20,181){\line(1,0){40}}
\put(21,155){$\bullet$}
\put(54,155){$\bullet$}
\put(15,178){$\bullet$}
\put(59,178){$\bullet$}
\put(36,198){$\otimes$}
\put(36,208){3}
\put(7,183){4}
\put(67,183){4}
\put(13,155){4}
\put(62,155){4}
\put(23,157){\line(-1,-1){12}}
\put(58,157){\line(1,-1){12}}
\put(19,180){\line(-1,0){16}}
\put(63,180){\line(1,0){16}}
\put(259,201){\line(-1,-1){21}}
\put(261,201){\line(1,-1){21}}
\put(244,157){\line(-1,4){6}}
\put(276,157){\line(1,4){6}}
\put(244,157){\line(1,0){33}}
\put(243,158){\line(2,5){17}}
\put(277,158){\line(-2,5){17}}
\put(239,180){\line(5,-3){36}}
\put(281,180){\line(-5,-3){36}}
\put(240,181){\line(1,0){40}}
\put(241,155){$\bullet$}
\put(274,155){$\bullet$}
\put(235,178){$\bullet$}
\put(279,178){$\bullet$}
\put(256,198){$\otimes$}
\put(256,208){3}
\put(227,183){4}
\put(287,183){4}
\put(233,155){4}
\put(282,155){4}
\put(243,157){\line(-1,-1){12}}
\put(278,157){\line(1,-1){12}}
\put(239,180){\line(-1,0){16}}
\put(283,180){\line(1,0){16}}
\end{picture}
\end{center}

\bigskip\noindent
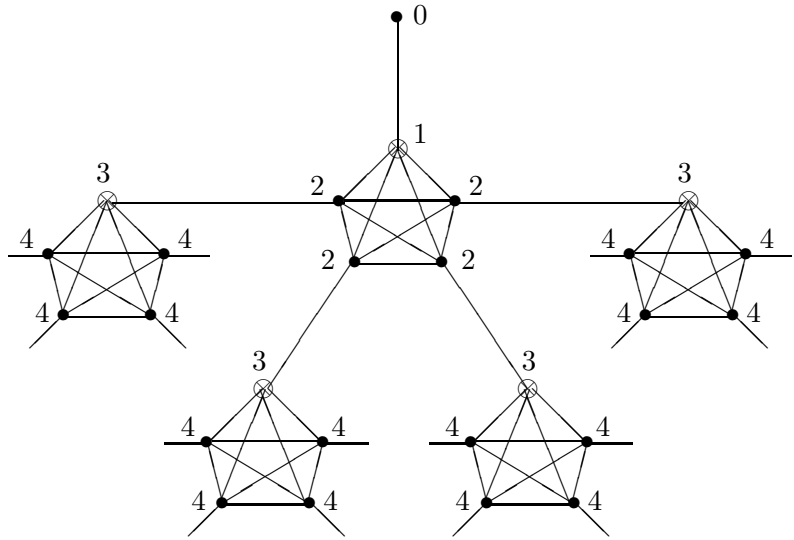
Figure 4. Comparison graph $\mathcal{K}_M$ when $M=5$. Numbers next
to the vertices indicate the corresponding states in the birth and death chain.

\clearpage

\begin{center}
\begin{picture}(300,270)
\put(150,110){\line(0,1){50}}
\put(25,110){\line(1,0){250}}
\put(25,160){\line(1,0){250}}
\put(50,160){\line(0,1){50}}
\put(100,160){\line(0,1){50}}
\put(200,160){\line(0,1){50}}
\put(250,160){\line(0,1){50}}
\put(50,110){\line(0,-1){50}}
\put(100,110){\line(0,-1){50}}
\put(200,110){\line(0,-1){50}}
\put(250,110){\line(0,-1){50}}
\put(25,185){\line(1,1){50}}
\put(75,185){\line(1,1){50}}
\put(175,185){\line(1,1){50}}
\put(225,185){\line(1,1){50}}
\put(25,35){\line(1,1){50}}
\put(75,35){\line(1,1){50}}
\put(175,35){\line(1,1){50}}
\put(225,35){\line(1,1){50}}
\put(34,194){\line(0,1){30}}\put(42,202){\line(0,1){30}}\put(58,218){\line(0,1){30}}\put(66,226){\line(0,1){30}}
\put(84,194){\line(0,1){30}} \put(92,202){\line(0,1){30}}\put(108,218){\line(0,1){30}}\put(116,226){\line(0,1){30}}
\put(184,194){\line(0,1){30}} \put(192,202){\line(0,1){30}} \put(208,218){\line(0,1){30}} \put(216,226){\line(0,1){30}}
\put(234,194){\line(0,1){30}} \put(242,202){\line(0,1){30}} \put(258,218){\line(0,1){30}} \put(266,226){\line(0,1){30}}
\put(34,44){\line(0,-1){30}} \put(42,52){\line(0,-1){30}} \put(58,68){\line(0,-1){30}} \put(66,76){\line(0,-1){30}}
\put(84,44){\line(0,-1){30}} \put(92,52){\line(0,-1){30}} \put(108,68){\line(0,-1){30}} \put(116,76){\line(0,-1){30}}
\put(184,44){\line(0,-1){30}} \put(192,52){\line(0,-1){30}} \put(208,68){\line(0,-1){30}} \put(216,76){\line(0,-1){30}}
\put(234,44){\line(0,-1){30}} \put(242,52){\line(0,-1){30}} \put(258,68){\line(0,-1){30}} \put(266,76){\line(0,-1){30}}
\put(147,108){$\bullet$}
\put(97,108){$\bullet$}
% top
\put(197,158){$\bullet$}
\put(147,158){$\bullet$}
\put(97,158){$\bullet$}
\put(73,160){*}
\put(223,160){*}
\put(193,180){*}
\put(97,208){$\bullet$}
\put(94,196){*}
\put(104,206){*}
% bottom left
\put(48,58){$\bullet$}
\put(42,54){*}
\put(52,62){*}
\put(93,80){*}
% bottom right
\put(197,58){$\bullet$}
\put(201,63){*}
\put(181,41){$\bullet$}
\put(174,35){*}
\put(178,20){*}
\end{picture}

\bigskip
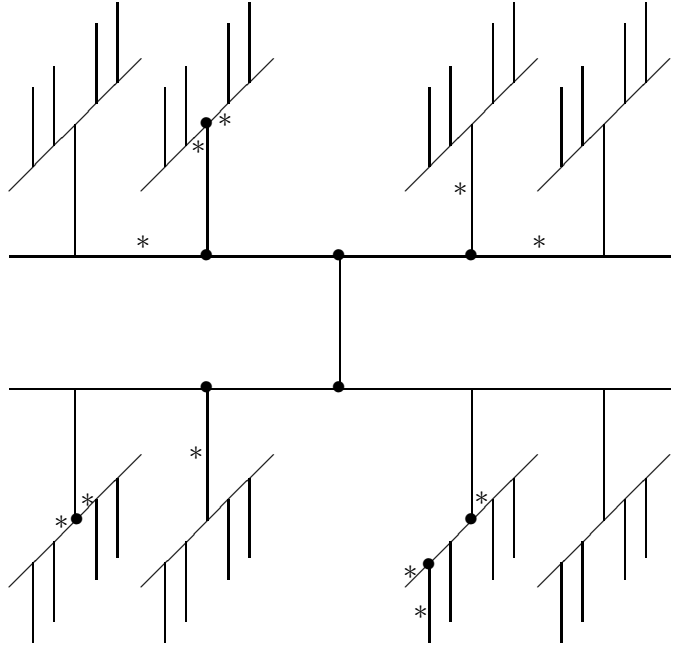
Figure 5. Finding exterior edges (*'s) on the big world graph, ${\cal B}$. Black dots indicate infected sites.
\end{center}

\clearpage

\begin{center}
\begin{picture}(200,200)
\put(20,100){\line(1,0){160}}
\put(40,120){\line(1,0){60}}
\put(60,140){\line(1,0){40}}
\put(80,160){\line(1,0){20}}
\put(100,80){\line(1,0){60}}
\put(100,60){\line(1,0){40}}
\put(100,40){\line(1,0){20}}
\put(39,100){\line(0,-1){20}}
\put(59,100){\line(0,-1){40}}
\put(79,100){\line(0,-1){60}}
\put(99,180){\line(0,-1){160}}
\put(119,100){\line(0,1){60}}
\put(139,100){\line(0,1){40}}
\put(159,100){\line(0,1){20}}
\put(16,98){$\bullet$}
\put(36,98){$\bullet$}
\put(56,98){$\bullet$}
\put(76,98){$\bullet$}
\put(96,98){$\bullet$}
\put(116,98){$\bullet$}
\put(136,98){$\bullet$}
\put(156,98){$\bullet$}
\put(176,98){$\bullet$}
\put(36,78){$\bullet$}
\put(56,78){$\bullet$}
\put(76,78){$\bullet$}
\put(96,78){$\bullet$}
\put(116,78){$\bullet$}
\put(136,78){$\bullet$}
\put(156,78){$\bullet$}
\put(36,118){$\bullet$}
\put(56,118){$\bullet$}
\put(76,118){$\bullet$}
\put(96,118){$\bullet$}
\put(116,118){$\bullet$}
\put(136,118){$\bullet$}
\put(156,118){$\bullet$}
\put(56,58){$\bullet$}
\put(76,58){$\bullet$}
\put(96,58){$\bullet$}
\put(116,58){$\bullet$}
\put(136,58){$\bullet$}
\put(56,138){$\bullet$}
\put(76,138){$\bullet$}
\put(96,138){$\bullet$}
\put(116,138){$\bullet$}
\put(136,138){$\bullet$}
\put(76,158){$\bullet$}
\put(96,158){$\bullet$}
\put(116,158){$\bullet$}
\put(76,38){$\bullet$}
\put(96,38){$\bullet$}
\put(116,38){$\bullet$}
\put(96,178){$\bullet$}
\put(96,18){$\bullet$}
\end{picture}
\end{center}
\bigskip\noindent
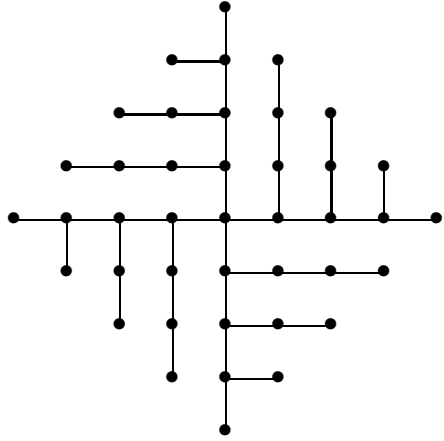
\centerline{Figure 6. Embedding a tree in $\mathbf{Z}^2$.}

\end{document}